\documentclass[12pt]{article}
\usepackage{graphicx}
\usepackage{here}

\newtheorem{thm}{Theorem}[section]
\newtheorem{lem}[thm]{Lemma}

\def\prf{\noindent {\em Proof. }}

\def\Qed{\hfill {\em Q.E.D.} \vskip .3in}

\def\PP{\vskip .1in \noindent}
\def\arg{{\rm arg~}}
\def\Re{{\rm Re~}}

\makeatletter
\def\theequation{\thesection.\@arabic\c@equation}
\makeatother

\begin{document}

\title {Zeros of Sections of the Binomial Expansion}
\author{
Svante Janson\footnote{Uppsala University, Sweden} \\
Timothy S. Norfolk\footnote{The University of Akron, U.S.A.}\\
}

\date{}
\maketitle

\noindent {\em Dedicated to Richard S. Varga, on the occasion of his 80th
birthday.}

\begin{abstract}
We examine the asymptotic behavior of the zeros of sections of the binomial
expansion. That is, we consider the distribution of
zeros of $B_{r,n}(z) = \sum_{k=0}^r {n \choose k} z^k$,
where $1 \le r \le n$.
\end{abstract}

\section{Preliminaries}
\setcounter{equation}{0}
\PP A problem of great interest in the classical Complex Function Theory is the
following:

\PP Given a function $f(z) = \sum_{k=0}^\infty a_k z^k$, analytic at $z=0$,
determine the asymptotic distribution of the zeros of the {\em partial sums}
$s_n(z) = \sum_{k=0}^n a_k z^k$.

\PP Some contributors to this area include Jentzsch \cite{J}, who explored the
problem for a finite radius of convergence; Szeg\H o \cite{Sz}, who explored
the exponential function $e^z$; Rosenbloom \cite{R}, who discussed the
angular distribution of zeros using potential theory, and applied his work to
sub-class of the confluent hypergeometric functions; Erd\H os and Tur\'an
\cite{ET}, who used minimization techniques to discuss angular distributions
of zeros; Newman and Rivlin \cite{NR1, NR2}, who related the work of Szeg\H o
to the Central Limit Theorem;
Edrei, Saff and Varga \cite{ESV}, who gave a thorough analysis
for the family of Mittag-Leffler functions; Carpenter, Varga and Waldvogel
\cite{CVW}, who refined the work of Szeg\H o; and Norfolk \cite{No1, No2},
who refined the
work of Rosenbloom on the confluent hypergeometric functions and a related
set of integral transforms.

\PP In this paper, we will analyze the behavior of the zeros of sections
of the binomial expansion, that is
\begin{equation}
B_{r,n}(z) = \sum_{k=0}^r {n \choose k} z^k~,1 \le r \le n~.
\end{equation}

\PP This investigation not only fits into the general theme of the works cited,
but also arises from matroid theory. Specifically (cf. \cite{W}), the
{\em univariate reliability polynomial} for the uniform matroid $U_{r,n}$
is given by
\begin{equation}
{\rm Rel}_{r,n} (q) = (1-q)^n B_{r,n} \left( \frac{q}{1-q} \right)
= \sum_{k=0}^r {n \choose k} q^k (1-q)^{n-k}~,
\end{equation}
which can be written as ${\rm Rel}_{r,n}(q) = (1-q)^{n-r} H_{r,n} (q)$,
where
\begin{equation}\label{rel}
H_{r,n}(q) = \sum_{k=0}^r {n \choose k} q^k (1-q)^{r-k}
= (1-q)^r B_{r,n} \left( \frac{q}{1-q} \right)~.
\end{equation}
\PP Some special cases are easy to analyze, and may thus be dispensed with. In
particular,
\begin{enumerate}
\item $B_{1,n} (z) = 1 + nz$, which has its only zero at $z =
-\frac{1}{n}$.
\item $B_{n,n} (z) = (1+z)^n$, which clearly has a zero of multiplicity $n$
at $z= -1$.
\item $B_{n-1,n} (z) = (1+z)^n - z^n$. Noting that this polynomial cannot have
positive zeros, we obtain the zeros $z = \frac{\omega^k}{1-\omega^k}$,
for $1 \le k \le n-1$, where $\omega = \exp \left( \frac{2\pi i}{n}
\right)$ is the principal $n$-th root of unity,
all of which lie on the vertical line ${\rm Re}~z 
= - \frac{1}{2}$.
\end{enumerate}

\PP In what follows, we will therefore focus on the cases
$1 \le r < n-1$, and give two collections of results. The first are concerned
with bounding regions for the zeros of $B_{r,n}(z)$, the rest with
convergence results.

\PP We note that this problem was investigated independently by Ostrovskii
\cite{Os}, who obtained many of the results that we present here.
The methods used there involved using a bilinear transformation to convert
the problem to an integral formulation. This choice of formulation makes the
proofs more involved and requires some additional constraints. By contrast,
we claim that our methods given here flow directly from the structure of the
problem, and yield additional results, in terms of additional bounds on the
zeros, and limiting cases. The paper \cite{Os} also gives a result on the
spacing of the zeros on the limit curve, using classical potential-theoretic
methods. We do not duplicate that result here, but give formulations in terms
of specific points on the curve.

\PP The methods used generate a set of constants and related
limit curves for $ 0 < \alpha < 1$,
defined by
\begin{equation}\label{cur1}
\frac{1}{2} \le K_\alpha = \alpha^\alpha (1-\alpha)^{1-\alpha} < 1~,
\end{equation}
\begin{equation}\label{cur2}
C_\alpha = \left\{ z~:~\frac{|z|^\alpha}{|1+z|} = K_\alpha
,~|z| \le \frac{\alpha}{1-\alpha} \right\} ~,
\end{equation}
and
\begin{equation}\label{cur3}
C_\alpha^\prime = \left\{ z~:~\frac{|z|^\alpha}{|1+z|} = K_\alpha
,~\frac{\alpha}{1-\alpha} \le |z| \right\}~.
\end{equation}
The properties of these curves are outlined in Lemma \ref{lem1}. Section 3
also presents bounds which are used to simplify the proofs of some of the
results presented here.

\section{Main Results}
\setcounter{equation}{0}

\PP As discussed above, we begin with a theorem on bounds of the zeros
of $B_{r,n}(z)$, and follow with results on convergence of those zeros.

\begin{thm}\label{thm1}
Let $r,n$ be positive integers, with $1 \le r < n-1$, and let
$z^*$ be any zero of $B_{r,n}(z) = \sum_{k=0}^r
{n \choose k} z^k$.

\PP Then, $z^*$ lies in a region defined by the intersection of two circles
and a
plane closed curve, to the right of a vertical line. Specifically,
\begin{equation}
|z^*| \le \frac{r}{n+1-r}~,
\end{equation}
\begin{equation}
\left| z^* - \frac{\gamma^2}{1-\gamma^2} \right|
\le \frac{\gamma}{(1-\gamma^2 )}~,
{\rm ~where~} \gamma = \frac{r}{n-1}~,
\end{equation}
\begin{equation}
{\rm Re~} z^* > -\frac{1}{2}~,
\end{equation}
and $z^*$ lies exterior to the curve $C_{r/n}$, as defined in (\ref{cur1},
\ref{cur2}).
\end{thm}

\prf We begin by considering the ratio of coefficients
\begin{equation}
\frac{{n \choose k}}{{n \choose k-1}} = \frac{n-k+1}{k}~,
\end{equation}
which is decreasing in $k$.

\PP Hence, writing $B_{r,n}\left( \frac{r}{n-r+1} z \right)
=\sum_{k=0}^r a_k z^k$, we have that
\begin{displaymath}
\frac{a_k}{a_{k-1}} = \frac{n-k+1}{k} \cdot \frac{r}{n-r+1} \ge 1~.
\end{displaymath}
That is, $\{ a_k\}_{k=0}^r$ is non-decreasing, so by the Enestr\"om-Kakeya
Theorem (\cite{He}, p. 462), the zeros of this polynomial satisfy $|z| \le 1$.
Hence, the zeros of $B_{r,n}(z)$ satisfy
$|z| \le \frac{r}{n-r+1}$.

\PP For the second bounding circle, we refer to Wagner \cite{W}, where it is
shown, again using the Enestr\"om-Kakeya Theorem, that the zeros of
$H_{r,n}(q)$ as given in (\ref{rel}), lie in the annulus
\begin{displaymath}
\frac{1}{n-r} \le |q| \le \frac{r}{n-1}~.
\end{displaymath}

\PP Since $z = -1$ is clearly not a zero of $B_{r,n} (z)$ for $r < n$,
we may make
the substitution $z = \frac{q}{1-q}$ (or equivalently
$q = \frac{z}{1+z}$) in (\ref{rel}), which shows
immediately that $H_{r,n}(q) = (1+z)^{-r} B_{r,n}(z)$, from which
one obtains
\begin{equation}\label{wag}
\left| \frac{z}{1+z} \right| \le \frac{r}{n-1} =: \gamma~.
\end{equation}

\PP Writing this last inequality in terms of the real and imaginary parts
of $z$ yields the claimed result.

\PP Noting that (\ref{wag}) implies that
$\left| \frac{z}{1+z} \right| < 1$,
yields the half-plane ${\rm Re~} z > -\frac{1}{2}$, as claimed.

\PP For the final bound, we mimic the analysis of Buckholtz \cite{Bu} on
the partial sums of $e^z$, and write
\begin{equation}\label{eq1.1}
(1+z)^{-n} B_{r,n}(z) = 1 - \frac{z^r}{(1+z)^n} \cdot
R_{r,n} (z)~,
\end{equation}
where
\begin{equation}\label{rem}
R_{r,n}(z) = \sum_{k=r+1}^n {n \choose k} z^{k-r}~
= z^{n-r} B_{n-r-1,n} \left( \frac{1}{z} \right).
\end{equation}

\PP For clarity, we set $\beta = r/n$. Inside and on the curve
$C_\beta$ (\ref{cur1},\ref{cur2}), we have
$|z| < \frac{\beta}{1-\beta }$
and $\left| \frac{z^r}{(1+z)^n} \right| \le K_\beta^n$,
where $K_\beta$ is defined in (\ref{cur1}). This, with the upper bound
of Lemma \ref{lem3} yields
\begin{equation}\label{eq1.2}
\left| (1+z)^{-n} B_{r,n}(z) \right| \ge 1-\left| \frac{z^r}{(1+z)^n} \right| \cdot
\left| R_{r,n}(z) \right| \\
 > 1 - K_\beta^n \cdot K_\beta^{-n} = 0 ~, \\
\end{equation}
which is the desired result. \Qed

\PP Note that
the second bounding circle occurring in this result,
namely
\begin{displaymath}
\left| z - \frac{\alpha^2}{1-\alpha^2} \right| = \frac{\alpha}{1-\alpha^2}~,
\end{displaymath}
intersects the negative real axis at
$z = -\frac{\alpha}{1+\alpha }$.
This circle is contained in the first, namely $|z| = \frac{\alpha}{1-\alpha}$, and
both meet at the common point $z = \frac{\alpha}{1-\alpha}$.

\begin{figure}[htp]
\includegraphics[angle=-90,width=4.3in]{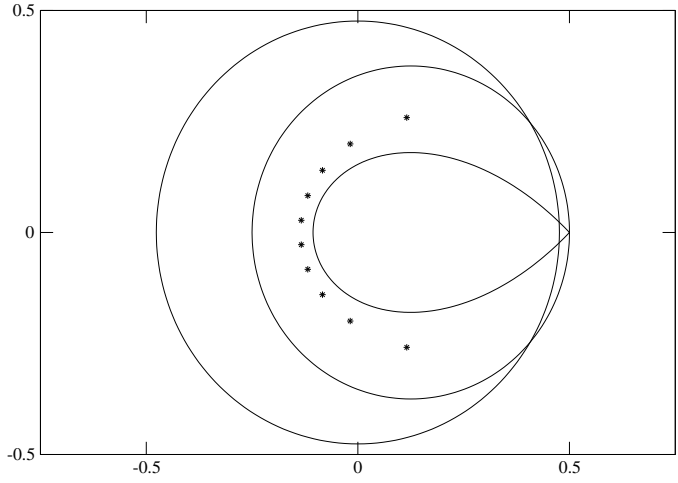}
\caption{The bounding curves and zeros for $r=10$, $n=30$}
\label{fig1}
\end{figure}

\PP The limiting case $|z| = \frac{\alpha}{1-\alpha}$ 
corresponding to the first
bounding circle, and the bounding half-plane ${\rm ~Re~} z > -\frac{1}{2}$ both
appear in \cite{Os}, with proofs that require significantly more detailed
derivation. The bounding curves and associated zeros for the case $r=10$
and $n = 30$ are illustrated in figure \ref{fig1}.

\PP We now use these results, and the bounds from the proof, to discuss 
some convergence results.

\begin{thm}\label{thm2}
Suppose that $1 \le r_j < n_j-1$ for all $j$, that
$\lim_{j \to \infty} n_j = \infty$, and that
\begin{displaymath}
\lim_{j \to \infty} \frac{r_j}{n_j} = \alpha,~0 < \alpha < 1~.
\end{displaymath}
Then
\begin{enumerate}
\item The zeros of $\{B_{r_j,n_j}(z)\}$ converge uniformly to points
of the curve $C_\alpha$, i.e.
\begin{displaymath}
\sup_{z:B_{r_j,n_j}(z)=0} d(z,C_\alpha ) \to 0~,
\end{displaymath}
where $d(z,C_\alpha ) = \inf_{\zeta \in C_\alpha} |z-\zeta |$ is the
distance from $z$ to $C_\alpha$,

and

\item Each point of $C_\alpha$ is a limit point of zeros of
$\left\{ B_{r_j,n_j} (z) \right\}_{j=1}^\infty$.
\end{enumerate}
\end{thm}

\prf Set $\beta_j = {r_j}/{n_j}$, so that
$\lim_{j \to \infty} \beta_j = \alpha$. Using (\ref{eq1.1}), the zeros
of $B_{r_j,n_j}(z)$
then satisfy
\begin{equation}\label{poly1}
\frac{z^{r_j}}{(1+z)^{n_j}} \cdot R_{r_j,n_j}(z) = 1~.
\end{equation}

\PP Using Theorem \ref{thm1}, Lemma \ref{lem1}
and Lemma \ref{lem3}, these zeros lie outside
the curve $C_{\beta_j}$, and thus satisfy
$\nu \beta_j < X_{\beta_j} \le |z| \le \frac{\beta_j}{1-\beta_j}$,
where $-X_{\beta_j}$ is the intersection of the curve $C_{\beta_j}$ with the
negative real axis, and $\nu$ is the unique positive solution to
$xe^{1+x} = 1$.

\PP Hence,
\begin{equation}\label{ineq}
\frac{\nu r_j}{n_j(r_j+1)} \le \frac{K_{\beta_j}^{n_j}
\left| R_{r_j,n_j}(z) \right|}
{\sum_{k=r_j+1}^{n_j} {n_j \choose k} {\beta_j}^k (1-\beta_j)^{n_j-k}} \le 1~,
\end{equation}
for this region. Note that the sum in the denominator above converges to $1/2$
by the Central Limit Theorem.

\PP Consequently,
$\lim_{j \to \infty} \left| R_{r_j,n_j}^{1/n_j}(z) \right| = K_\alpha^{-1}$
uniformly on the
set in question. Taking
moduli and $n_j$-th roots in (\ref{poly1}), we observe that
the zeros of $B_{r_j,n_j}(z)$ must satisfy
\begin{equation}\label{asym}
\frac{|z|^{\beta_j}}{|1+z|} |R_{r_j,n_j}(z)|^{1/n_j} = 1~.
\end{equation}

\PP Since $\beta_j \to \alpha$, this establishes that every limit point of
a sequence of zeros of $B_{r_j,n_j}(z)$ lies on $C_\alpha$. Since, by Theorem
\ref{thm1}, the zeros lie in a compact set, it follows that the zeros converge
uniformly to points of $C_\alpha$.

\PP For the second claim, fix any $\zeta \in C_\alpha$ with
$\zeta \ne z_\alpha = \alpha/(1-\alpha )$. Then $|\zeta | < z_\alpha$, so we
may take a small neighborhood $D$ of $\zeta$ such that $0 < |z| < z_\alpha$
for $z \in {\overline D}$. Consequently, for $j$ sufficiently large,
$|z| < z_{\beta_j}$ for all $z \in {\overline D}$, and it follows from Lemma
\ref{lem3} and the Central Limit Theorem, that
\begin{displaymath}
\left| R_{r_j, n_j}(z) \right|^{1/n_j} \to K_\alpha^{-1}~,
\end{displaymath}
uniformly on ${\overline D}$.

\PP In particular, for large $j$, $R_{r_j,n_j}(z) \ne 0$ on ${\overline D}$,
so we may fix an analytic branch of $R_{r_j,n_j}^{1/n_j}(z)$ in $D$. Letting
$\theta_j = \arg (R_{r_j,n_j}^{1/n_j}(z))$ (with arguments in the
range $(0,2\pi )$), we then have
\begin{displaymath}
e^{-i \theta_j} R_{r_j,n_j}(z) \to K_\alpha^{-1}~,
\end{displaymath}
uniformly on compact subsets of $D$.

\PP By shrinking $D$, we may assume that the latter limit holds uniformly on
$D$.  Furthermore, we may assume that $0 < \arg (z) < 2 \pi$ for $z \in D$,
and thus the powers $z^{\beta_j}$ and $z^\alpha$ are well-defined in $D$. Hence,
\begin{equation}\label{svan}
\frac{z^{\beta_j}}{1+z} R_{r_j,n_j}^{1/n_j}(z) - \frac{z^\alpha}{1+z}
K_\alpha^{-1} e^{i \theta_j} \to 0~,
\end{equation}
uniformly on $D$.

\PP Since the mapping $w = \frac{z^\alpha}{1+z} K_\alpha^{-1}$ maps $C_\alpha$
onto an arc of the unit circle, it maps $D \cap C_\alpha$ onto a subarc. Thus,
for $j$ sufficiently large, there exists an integer $p_j$ such that
$\frac{z^\alpha}{1+z} K_\alpha^{-1} e^{i \theta_j} = e^{2\pi i p_j /n_j}$
for some $z = \zeta_j \in D \cap C_\alpha$. We may further assume that
$\zeta_j \to \zeta$. It now follows from Hurwitz' theorem and (\ref{svan}) that,
for $j$ sufficiently large,
\begin{displaymath}
\frac{z^{\beta_j}}{1+z} R_{r_j,n_j}^{1/n_j}(z) - e^{2 \pi i p_j / n_j}
\end{displaymath}
has a zero $z_j \in D$. Each such zero $z_j$ satisfies (\ref{poly1}), and so
by (\ref{eq1.1}), is a zero of $B_{r_j,n_j}(z)$. This proves that every point on
$C_\alpha$ is a limit point of zeros of $\left\{ B_{r_j,n_j}(z) \right\}$.

\Qed

\PP We note that, thanks to (\ref{rem}), the non-trivial zeros of
$R_{r_j,n_j}(z)$
converge uniformly to all points which lie on the curve $C_\alpha^\prime$,
as defined in (\ref{cur3}).

\PP This result also appears in \cite{Os}, using more elaborate asymptotics.
The analysis presented requires a deletion of a neighborhood of the 
singular point
$z_\alpha = \frac{\alpha}{1-\alpha }$. Consideration of the results of
Lemma \ref{lem3} shows that this is not necessary with our methods.

\PP The remaining results presented here do not appear in the literature.

\PP The asymptotic expansions in the proof of Theorem \ref{thm2}
immediately give
the following result on the rate of convergence. We note that, as shown
in \cite{CVW} in the case of the exponential function, this rate is
best possible.

\begin{thm}\label{thm3}
Fix $0 < \delta < 1$. Then, there exists a constant $c$, depending only on $\delta$,
such that, if
$r,n$ are large, and $0 < \delta < \frac{r}{n} < 1-\delta$,
for any zero $z^*$ of $B_{r,n}(z)$
\begin{displaymath}
\min_{\zeta \in C_{r/n}} |z^*-\zeta| \le \frac{c}{|z^*-\frac{r}{n-r}|}\cdot \frac{\ln n}{n}~.
\end{displaymath}

\PP Additionally, proximity to the singular point
$z_{r/n} = \frac{r}{n-r}$ is of order $O \left( \frac{1}{\sqrt{n}}\right)$.
\end{thm}

\prf Set $\beta = r/n$ From (\ref{ineq}), we obtain the approximation
\begin{equation}\label{svan2}
\left| R_{r,n}^{1/n} (z) \right| \cdot K_\beta = 1+G_{r,n}(z)
\cdot \frac{\ln n}{n}~,
\end{equation}
where $G_{r,n}(z)$ is uniformly bounded in a region containing the zeros.

\PP Let $z^*$ be a zero of $B_{r,n}(z)$, and let $\zeta$ be the point on
$C_\beta$ closest to $z^*$. Note that $|\zeta-z^*| = o(1)$ as a consequence of
Theorem \ref{thm2}, as applied to sequences for which $\beta$ converges.
Note that the curve $C_\beta$ is asymptotically a pair of straight lines
at angle $\pi/4$ to the real axis close to the point
$z_\beta = \beta /(1-\beta )$.
Hence, if $z^*$ is close to $z_\beta$, by Theorem \ref{thm1}, it must lie
in the wedges between these lines and the vertical line $\Re z = z_\beta$,
from which $|z^*-z_\beta | = O(|\zeta - z_\beta |)$.

\PP Note that $z^*$ satisfies (\ref{asym}), without the subscript $j$, and thus,
by (\ref{svan2}), we have
\begin{displaymath}
\frac{|z^*|}{|1+z^*|}\cdot K_\beta^{-1} = \left( 1 + G_{r,n}(z)
\cdot \frac{\ln n}{n} \right)^{-1}~.
\end{displaymath}
Expanding $F(z) = \ln (K_\beta^{-1} |z|^\beta/|1+z|)
= \Re \ln(K_\beta^{-1} z^\beta / (1+z))$ as a Taylor series centred at $\zeta$
(noting that $F(\zeta) = 0$), we find that
\begin{displaymath}
|z^*-\zeta | = O \left( \left| \frac{\zeta (1+\zeta)}{\beta - (1-\beta ) \zeta}
\cdot G_{r,n}(z)\cdot \frac{\ln n}{n} \right| \right)
= O \left( \frac{1}{|z_\beta - \zeta|}\cdot \frac{\ln n}{n} \right)~.
\end{displaymath}

This not only gives the desired result, but shows that,
as expected, the rate of convergence is worst for those points closest to
the singular point
$z_\beta = \frac{\beta}{1-\beta }$.

\PP To discuss the convergence at the singular point, we take an approach
similar to that used
for the exponential function in \cite{NR1,NR2} and for the Mittag-Leffler
functions in \cite{ESV}.
For convenience, we set
$\mu = n \beta = r$, and $\sigma^2 = n \beta (1-\beta)$. Then,
\begin{displaymath}
f_{r,n}(w) =
(1-\beta)^n B_{r,n} \left( \frac{\beta e^{w/\sigma}}{1-\beta} \right)
= \sum_{k=0}^r {n \choose k} \beta^k (1-\beta )^{n-k} e^{kw/\sigma}~,
\end{displaymath}
which is a truncated moment generating function for a binomial distribution
with mean $\mu$ and variance $\sigma$. Using the Central Limit Theorem,
\begin{displaymath}
f_{r,n}(w) \approx \frac{1}{\sqrt{2\pi}\sigma} \int_{-\infty}^\mu
e^{-\frac{1}{2} \left( \frac{t-\mu}{\sigma} \right) + \frac{tw}{\sigma}} dt~.
\end{displaymath}
Making the substitution $s = \frac{t-\mu-\sigma w}{\sqrt{2}\sigma}$ yields
\begin{displaymath}
e^{-\mu w/\sigma - w^2/2} f_{r,n}(w) \approx
\frac{1}{\sqrt{\pi}} \int_{-\infty}^{-w/\sqrt{2}} e^{-s^2} ds
= \frac{1}{2} {\rm erfc} \left( \frac{w}{\sqrt{2}} \right)~,
\end{displaymath}
the complementary error function. Thus, given the zero $\chi$ of ${\rm erfc}(z)$
which is closest to the origin, there must exist a zero $z^*$ of $B_{r,n}(z)$
for which
\begin{displaymath}
z^* \approx \frac{\beta e^{\sqrt{2}\chi / \sigma}}{1-\beta}
\approx \frac{\beta}{1-\beta} + \sqrt{\frac{2\beta}{(1-\beta)^3}}\cdot
\frac{\chi}{\sqrt{n}}~,
\end{displaymath}
the desired result.
\Qed

\begin{figure}[hbp]
\includegraphics[angle=-90,width=4.3in]{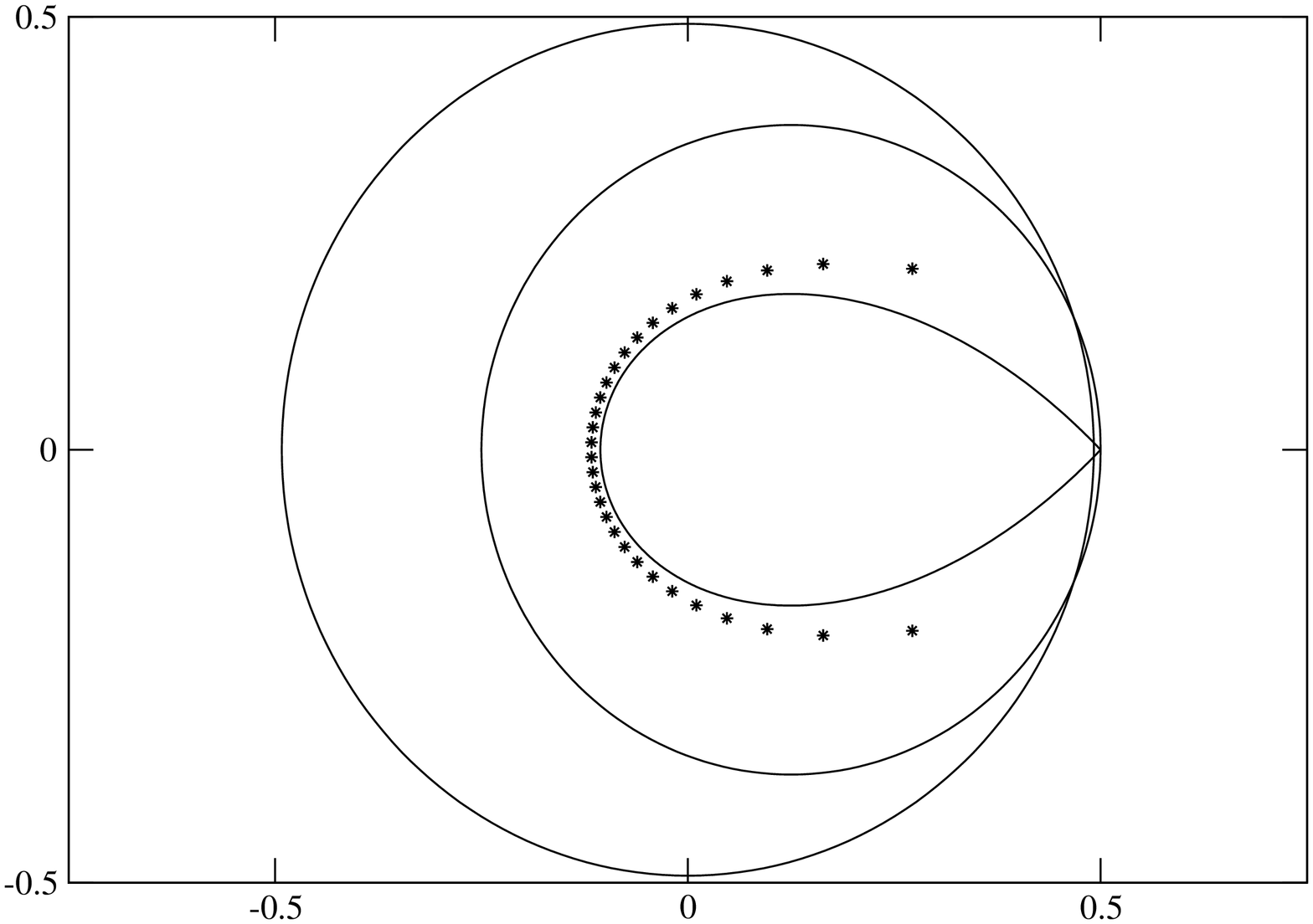}
\caption{The bounding curves and zeros for $r=30$, $n=90$}
\label{fig2}
\end{figure}

\PP The figures \ref{fig1} and \ref{fig2} show the zeros, bounding curve
and bounding circles for the cases $r=10, n=30$ and $r=30, n=90$
respectively. Since the ratio $r/n$ is the same in both cases, they serve
to illustrate both the rate of convergence of the zeros to the limit curve,
and the rate of convergence of the bounding circles.

\begin{figure}[htp]
\includegraphics[angle=-90,width=4.3in]{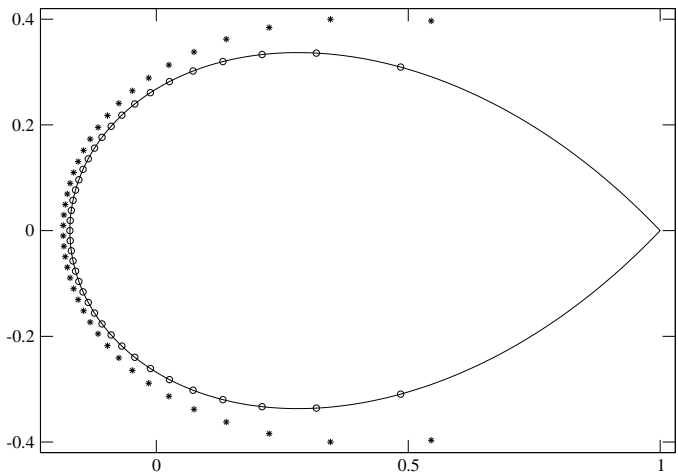}
\caption{The curve $C_{1/2}$, points $\{\zeta_{p,80}\}$
and zeros for $r=40$, $n=80$}
\label{fig3}
\end{figure}

\PP Figure \ref{fig3} shows the zeros for the case $r=40, n=80$, as well as the
curve $C_{1/2}$ and the approximating points on the curve.

\PP It should be noted at this point that, due to the structure of the
coefficients of these polynomials, direct computation of the zeros
for significantly higher degrees suffers due to numerical instability.

\PP We conclude by considering the limiting cases $\alpha = 0$ and
$\alpha = 1$. The trivial
result for $\alpha = 0$, given the radius $\frac{r}{n+1-r}$ of the bounding
circle, is that
all zeros converge uniformly to $0$ in this case. However, a slight modification
gives a much more interesting result.

\begin{thm}\label{thm4}
Suppose that $\lim_{j \to \infty} r_j = \infty$ and that
$\lim_{j \to \infty} \frac{r_j}{n_j} = 0$.

\noindent Then, the limit points of the zeros of $\{ B_{r_j,n_j}
( \frac{r_j z}{n_j-r_j} ) \}_{j=1}^\infty$
are precisely the points of the
Szeg\H o curve $|ze^{1-z}| = 1$, $|z| \le 1$.
\end{thm}

\prf 
With the given normalization, the results of Theorem \ref{thm1}
yield that the zeros of the normalized polynomial above satisfy
\begin{equation}\label{eq4.1}
1 = \left( \frac{r_j}{n_j-r_j} \right)^{r_j} K_{r_j/n_j}^{-n_j}
\frac{z^{r_j}}{\left( 1+\frac{r_jz}{n_j-r_j} \right)^{n_j}} h(z)
{\rm ~and~} |z| \le 1~,
\end{equation}
where

\begin{equation}\label{eq4.2}
h(z) =
\sum_{k=r_j+1}^{n_j} {n_j \choose k} \left( \frac{r_j}{n_j} \right)^k \left(
1- \frac{r_j}{n_j} \right)^{n_j-k_j} z^{k-r_j}.
\end{equation}

\PP Noting that
\begin{displaymath}
\left( \frac{r_j}{n_j-r_j} \right)^{r_j} K_{r_j/n_j}^{-n_j}
= \left( 1 - \frac{r_j}{n_j} \right)^{-n_j}~,
\end{displaymath}
we may use standard expansions to convert (\ref{eq4.1}) to the form
\begin{equation}\label{eq4.3}
1 = (ze^{1-z+g(z)})^{r_j} h(z)~,
\end{equation}
where $|g(z)| \le \frac{3r}{n}$ uniformly in the unit disk.

\PP Considering points inside and on the curve $|ze^{1-z}| = e^{-3r_j/n_j}$,
and noting that $|h(z)| \le h(1) < 1$ on the unit disk,
we may repeat the analysis of (\ref{eq1.2}) to deduce that the zeros are
uniformly bounded away from zero by $|z| \ge \eta > 0$.
This implies that we may repeat the bounding process of
Lemma \ref{lem3} to deduce that $h^{1/r_j}(z) \to 1$ uniformly
in $\eta \le |z| \le 1$, defining the roots by a cut along the positive
real axis.
This establishes the desired result.
\Qed

\PP Finally, we consider the other limiting case.

\begin{thm}\label{thm5}
Suppose that $\lim_{j \to \infty} r_j = \infty$ and
$\lim_{j \to \infty} \frac{r_j}{n_j} = 1$.

\noindent Then, the limit points of the zeros of the
polynomials $\left\{ B_{r_j,n_j}(z) \right\}_{j=1}^\infty$
are precisely the points of the line ${\rm Re~}z = -\frac{1}{2}$.
\end{thm}

\prf 
As in the previous proofs, we write the
equation for the zeros as
\begin{displaymath}
1  = \frac{z^{r_j}}{(1+z)^{n_j}} R_{r_j,n_j}(z)~.
\end{displaymath}

\PP We again use the bounds of Lemma \ref{lem3} and obtain the desired
result, using the fact that $\lim_{\alpha \to 1^-} K_\alpha = 1$.
\Qed

\section{Technical Results}
\setcounter{equation}{0}

\PP Here we give the properties and inequalities necessary for the main results,
beginning with the properties of the bounding curves.

\begin{lem}\label{lem1}
Fix $0 < \alpha < 1$, and let
\begin{equation}
K_\alpha = \alpha^\alpha (1-\alpha )^{1-\alpha}
\end{equation}
and
\begin{equation}
C_\alpha = \left\{ z~:~\frac{|z|^\alpha}{|1+z|} = K_\alpha ,~
|z| \le \frac{\alpha}{1-\alpha} \right\}~.
\end{equation}
Then,
\begin{enumerate}
\item $\frac{1}{2} \le K_\alpha < 1$, $\lim_{\alpha \to 0^+}
K_\alpha = 1$, $\lim_{\alpha \to 1^-} K_\alpha = 1$.
\item $C_\alpha$ is a simple, smooth closed curve, symmetric with respect
to the real axis, starlike with respect to $z=0$, which passes through
$z_\alpha = \frac{\alpha}{1-\alpha }$.
\item The intersection of $C_\alpha$ with the negative real axis
occurs at $z = -X_\alpha$,where
$\nu \alpha < X_\alpha < \frac{1}{2}$ and $\nu = 0.278 \cdots$
is the unique positive root of
$xe^{1+x} = 1$.
\item $X_\alpha \le |z|$ and $|z| \le z_\alpha$ for any $z \in C_\alpha$,
with the latter equality holding only at $z = z_\alpha$.
\end{enumerate}
\end{lem}

\prf \begin{enumerate} \item A simple calculation gives the limits. Taking
derivatives yields
\begin{displaymath}
\frac{dK_\alpha}{d\alpha} = K_\alpha \ln \left( \frac{\alpha}{1-\alpha}
\right)~,
\end{displaymath}
which shows that $K_\alpha$ is decreasing on $\left( 0, \frac{1}{2}
\right)$ and increasing on $\left( \frac{1}{2}, 1 \right)$. Calculating
$K_{1/2}$ directly gives the equality.

\item Clearly, the definition shows that $C_\alpha$ is closed and symmetric,
and direct calculation shows that it passes through the point
$z_\alpha = \alpha/(1-\alpha )$.

\PP We write $z = re^{i\theta}$, and set
\begin{equation}
c_\theta (r) = \frac{|z|^\alpha}{|1+z|} =
\frac{r^\alpha}{\sqrt{1+2r\cos \theta + r^2}}~.
\end{equation}
Clearly, $c_\theta (0) = 0$ and $\lim_{r \to \infty} c_\theta (r) = 0$.

\PP For $\theta = 0$, we have
\begin{displaymath}
c_0^\prime (r) = \frac{r^{\alpha-1}}{(1+r)^2} [ \alpha - (1-\alpha )r ]~,
\end{displaymath}
which shows that the given point is the only positive real
value satisfying the
equation.

\PP For $0 < \theta < \pi$, we have
\begin{displaymath}
c_\theta^\prime (r) = r^{\alpha-1}{(1+2r \cos \theta +r^2)^{-3/2}}
[(\alpha-1) r^2 + (2\alpha -1) r \cos \theta + \alpha]~.
\end{displaymath}

\PP Since $\alpha - 1 < 0$, this derivative has exactly one positive root,
which is a maximum of the function. Further, a simple calculation shows that
\begin{displaymath}
c_\theta \left( \frac{\alpha}{1-\alpha} \right) > K_\alpha~,
\end{displaymath}
from which each such ray yields exactly one point on the curve, inside the
bounding circle, $|z| = \frac{\alpha}{1-\alpha}$.
Considering the defining function, this value of $r$
is clearly decreasing in $0 \le \theta < \pi$.
Hence, the curve is simple and starlike with respect to 0.

\PP Finally, for $\theta = \pi$, we have that
\begin{displaymath}
c_\pi^\prime (r) = \frac{r^{\alpha-1}}{(1-r)^2} [\alpha + (1-\alpha )r] > 0
\end{displaymath}
for $0 < r < 1$, and $\lim_{r \to 1^-} c_\pi (r) = \infty$, which gives
exactly one solution in this range.

\PP That these points are the only solutions within the bounding circle
can be deduced from the fact that $z \in C_\alpha$
if and only if
$\frac{1}{z} \in C_{1-\alpha}^\prime$.

\PP Examining the function $w = K_\alpha^{-1}\frac{z^\alpha}{1+z}$ using arguments
in the range $(0, 2\pi )$ shows that $C_\alpha$ maps onto the approriate
arc of the unit circle in the $w$-plane. This mapping is also one-to-one
along the arc $0 < \arg w < 2 \pi \alpha$, since $w^\prime \ne 0$ on the cut plane.
This fact is implicitly
used in the calculation of the rate of convergence.

\item The solution on the negative real axis is $-t = -X_\alpha$, and satisfies
\begin{displaymath}
\frac{t^\alpha}{1-t} = K_\alpha ~,
\end{displaymath}
which we write as
\begin{equation}
f(t) = t^\alpha + \alpha^\alpha (1-\alpha )^{1-\alpha} (t-1) = 0~.
\end{equation}
\PP Now, $f(t)$ is increasing, with $f(0) < 0$,
$f (X_\alpha ) = 0$, and
\begin{displaymath}
f \left( \frac{1}{2} \right) = \left( \frac{1}{2} \right)^\alpha - \frac{1}{2} K_\alpha
> \frac{1}{2} (1-K_\alpha ) > 0~,
\end{displaymath}
from which $X_\alpha < \frac{1}{2}$ follows immediately.

\PP To show that $\nu \alpha < X_\alpha$, we consider
\begin{equation}\label{eq3.1}
f(\nu \alpha ) = \alpha^\alpha ( \nu^\alpha - (1- \nu \alpha )
(1-\alpha )^{1-\alpha} )~.
\end{equation}
and set
\begin{equation}
g(\alpha ) = \ln ((1- \nu \alpha )(1-\alpha )^{1-\alpha })~,
\end{equation}
which satisfies $g (0) = 0$, $g^\prime (0) = -\nu - 1$ and
\begin{equation}
g^{\prime\prime} (\alpha ) = \frac{ (\nu \alpha )^2 +
(\nu-2)(\nu \alpha ) + 1 - \nu^2 }
{(1-\alpha ) (1-\nu \alpha )^2} > 0~.
\end{equation}
The last inequality follows since the quadratic in the numerator has
discriminant $\nu^3(5\nu-4) < 0$, from Lemma \ref{lem1}, and so has no
real zeros.

\PP Hence, 
\begin{displaymath}
e^{g(\alpha)} > e^{-(\nu + 1) \alpha} = e^{\alpha \ln \nu } = \nu^\alpha~,
\end{displaymath}
and thus, by (\ref{eq3.1}), $f(\nu \alpha ) < 0$ for $0 < \alpha < 1$,
as desired. \Qed

\end{enumerate}

\PP We continue with a lemma required for one of the bounds.

\begin{lem}\label{lem2}
Let $f(z) = \sum_{k=0}^\infty b_k z^k$ satisfy
\begin{equation}
b_0 > b_1 \ge 0,~b_k \ge 0,~b_1b_{k-1}-b_0b_k \ge 0 {\rm ~for~} k \ge 1~.
\end{equation}
Then, $|f(z)| \ge \frac{b_0-b_1}{b_0+b_1} f(1)$ for $|z| \le 1$.
\end{lem}

\prf The conditions given imply that $\{ b_k \}$ is strictly decreasing, unless
$b_k = 0$ for $k \ge K$.
Let $r = \frac{b_1}{b_0} < 1$. Then, the conditions given show that
$b_k \le r b_{k-1}$ for $k \ge 1$. Hence, $f(z)$ is analytic for $|z| < \frac{1}{r}$,
and in particular in the closed unit disk.
Applying the Enestr\"om-Kakaya Theorem to the partial sums
$p_n(z) = \sum_{k=0}^n b_k z^k$ shows that all have their zeros
in the region $|z| > 1$, hence, by Hurwitz' Theorem, $f(z)$ cannot have any zeros
inside the unit disk. Thus, applying the Minimum Modulus Theorem, the minimum value
of $|f(z)|$ for $|z| \le 1$ must occur on the boundary.

\PP For $|z| = 1$, we have
\begin{equation}
\begin{array}{lcl}
|(b_0-b_1z) f(z)| & = & \left| b_0^2 + \sum_{k=1}^\infty (b_0b_k-b_1b_{k-1}) z^k \right| \\
& \ge & b_0^2 - \sum_{k=1}^\infty \left| (b_1b_{k-1}-b_0b_k) \right| \\
& = & b_0^2 - \sum_{k=1}^\infty b_1b_{k-1} + \sum_{k=1}^\infty b_0b_k \\
& = & b_0^2-b_1f(1)+b_0 (f(1)-b_0) \\
& = & (b_0-b_1) f(1) ~. \\
\end{array}
\end{equation}
Hence, we have
\begin{equation}
|f(z)| \ge \frac{(b_0-b_1)f(1)}{|b_0-b_1z|} \ge \frac{(b_0-b_1)f(1)}{b_0+b_1}~,
\end{equation}
the desired result. \Qed

\PP Finally, we have the estimates of the remainder term.

\begin{lem}\label{lem3} Given integers $1 \le r < n$, we set
$\beta = \frac{r}{n}$,
and consider the remainder term
\begin{equation}
R_{r,n}(z) = \sum_{k=r+1}^n {n \choose k} z^{k-r}~.
\end{equation}
Then, for $|z| \le \frac{\beta}{1-\beta}$, we have
\begin{equation}
\left| R_{r,n}(z) \right| \le K_\beta^{-n} \sum_{k=r+1}^n {n \choose k}
\beta^k (1-\beta )^{n-k} \le K_\beta^{-n}
\end{equation}
and
\begin{equation}
\left|  R_{r,n}(z) \right| \ge
\frac{|z|}{r+1} K_\beta^{-n}
\sum_{k=r+1}^n {n \choose k} \beta^k (1-\beta )^{n-k}~.
\end{equation}
\end{lem}

\prf Given that all coefficients are positive, we use the value of $K_\beta$
from (\ref{cur1}) and the bound on $|z|$ to deduce that
\begin{displaymath}
\left| R_{r,n}(z) \right| \le R_{r,n} \left( \frac{\beta}{1-\beta} \right)
= K_\beta^{-n} \sum_{k=r+1}^n {n \choose k} \beta^k (1-\beta )^{n-k}~.
\end{displaymath}

\PP The latter sum is clearly bounded by 1, using the binomial expansion. In fact,
using the Central Limit Theorem, it is asymptotically $1/2$ for $r$ and $n-r$
both large.

\PP For the lower bound, we consider
\begin{displaymath}
g(z) = \left( \frac{1-\beta }{\beta z} \right) R_{r,n} \left( \frac{\beta z}{1-\beta }
\right) = \sum_{k=0}^{n-r-1} b_k z^k~,
\end{displaymath}
where
\begin{displaymath}
b_k = {n \choose k+r+1} \left( \frac{\beta}{1-\beta } \right)^k~.
\end{displaymath}

\PP It is simple to show that $g(z)$ satisfies the conditions of Lemma \ref{lem2},
that
\begin{displaymath}
\frac{b_0-b_1}{b_0+b_1} = \frac{2n-r}{2r(n-r)+(2n-3r)} \ge \frac{1}{r+1}~,
\end{displaymath}
and finally that
\begin{displaymath}
g(1) = \left( \frac{1-\beta}{\beta} \right) K_\beta^{-n} \sum_{k=r+1}^n
{n \choose k} \beta^k (1-\beta )^{n-k}~.
\end{displaymath}

\PP Rewriting $R_{r,n}(z)$ in terms of $g(z)$ yields the result.
\Qed

\PP We would like to acknowledge Professor Alan Sokal, of New York University,
who suggested this problem in 2001, and independently deduced the form
of the limit curves.
\nocite{*}
\bibliographystyle{plain}
\bibliography{bin}

\end{document}